\documentclass[runningheads]{llncs}
\pdfoutput=1
\usepackage{graphicx}
\usepackage{amsmath}
\usepackage{amssymb}
\usepackage{todonotes}
\usepackage{bbm} % for the mathbb 1
\usepackage{longtable}
\usepackage{aliascnt} % Used for references like "theorem 3.4"
\usepackage{hyperref}
\usepackage{multirow} %for table
\DeclareMathOperator{\STAB}{STAB}
\DeclareMathOperator{\COL}{COL}
\DeclareMathOperator{\CUT}{CUT}
\DeclareMathOperator{\conv}{conv}
\DeclareMathOperator{\trace}{trace}
\DeclareMathOperator{\diag}{diag}

\newcommand{\COLE}{\COL^{\varepsilon}}
\newcommand{\eshmc}[1]{z_{mc}^{#1}}
\newcommand{\eshss}[1]{z_{ss}^{#1}}
\newcommand{\eshc}[1]{z_{c}^{#1}}
\newcommand{\eshcE}[1]{z_{c\varepsilon}^{#1}}
\newcommand{\esbmc}[1]{s_{mc}^{#1}}
\newcommand{\esbss}[1]{s_{ss}^{#1}}
\newcommand{\esbc}[1]{s_{c}^{#1}}

\newcommand{\XI}{X_{I}}
\newcommand{\GI}{G_{I}}
\newcommand{\bI}{b_{I}}
\newcommand{\kI}{k_{I}}
\newcommand{\yI}{y_{I}}
\newcommand{\yJ}{y}
\newcommand{\lambdaI}{\lambda_{I}}
\newcommand{\lambdaJ}{\lambda}

\newcommand{\opvt}{ {\transposedVec} }

\newcommand{\sumbI}{b}
\newcommand{\tI}{t_{I}}
\newcommand{\setX}{{\mathcal X}}
\newcommand{\setXmc}{{\mathcal X}^{E}}
\newcommand{\setXss}{{\mathcal X}^{S}}
\newcommand{\setXcol}{{\mathcal X}^{C}}

\newcommand{\R}{{\mathbb R}}
\newcommand{\Sym}[1]{{\mathcal S}_{#1}}
\newcommand{\Lagrangian}{{\mathcal L}}
\newcommand{\simplex}{{\Delta}}
\newcommand{\allones}[1]{\mathbbm{1}_{#1}}
\newcommand{\transposedVec}{T}
\newcommand{\transposedOp}{\top}
\newcommand{\norm}[1]{\left\lVert #1 \right\rVert}

\newcommand{\opA}{{\mathcal A}_I}
\newcommand{\opAt}{{\mathcal A}_I^{\transposedOp}}
\newcommand{\opAmc}{{\mathcal A}_I}
\newcommand{\opAmct}{{\mathcal A}_I^{\transposedOp}} %\leftidx{^{C}}
\newcommand{\opM}{{\mathcal M}_I}
\newcommand{\opMt}{{\mathcal M}_I^{\transposedOp}}
\newcommand{\opMmc}{{\mathcal M}_I}
\newcommand{\opMmct}{{\mathcal M}_I^{\transposedOp}} %\leftidx{^{C}}
\newcommand{\opP}{{\mathcal P}_I}
\newcommand{\opPt}{{\mathcal P}_I^{\transposedOp}}
\newcommand{\opD}{{\mathcal D}_I}

\newcommand{\hBdlExplanation}[1]{h_{#1}}
\newcommand{\yBdlExplanation}[1]{y_{#1}}
\newcommand{\gBdlExplanation}[1]{g_{#1}}
\newcommand{\eBdlExplanation}[1]{e_{#1}}
\newcommand{\XBdlExplanation}[1]{X_{#1}}

\newcommand{\yCenter}{\overline{y}}
\newcommand{\yTrial}{\widetilde{y}}

\newcommand{\hCenter}{\overline{h}}
\newcommand{\wTrial}{\widetilde{w}}
\newcommand{\vTrial}{\widetilde{v}}
\newcommand{\vI}{v_{I}}
\newcommand{\bdlSize}{r}
\newcommand{\bdlPar}{\mu}
\newcommand{\nrESC}{q}
\newcommand{\Bundle}{\mathcal{B}}

\begin{document}
\title{A Bundle Approach for SDPs with Exact Subgraph Constraints\thanks{This project has received funding from the European Union’s Horizon~2020 research and innovation programme under the Marie Sk\l{}odowska-Curie grant agreement No~764759 and the Austrian Science Fund (FWF): I 3199-N31 and P 28008-N35. We thank three anonymous referees for their constructive comments which substantially helped to improve the presentation of our material.}}

%The first author acknowledges support by the Austrian Science Fund (FWF): I 3199-N31.
\author{Elisabeth Gaar\inst{1}\orcidID{0000-0002-1643-6066} \and
Franz Rendl\inst{1}\orcidID{0000-0003-1578-9414}}

\authorrunning{E. Gaar and F. Rendl}
\institute{Alpen-Adria-Universit\"{a}t Klagenfurt, Institut f\"{u}r Mathematik, \\
    Universit\"{a}tsstr. 65-67, 9020 Klagenfurt, Austria \\
    \email{\{elisabeth.gaar,franz.rendl\}@aau.at}}
\maketitle              % typeset the header of the contribution
\begin{abstract}
The 'exact subgraph' approach was recently introduced as a hierarchical
scheme to get increasingly tight semidefinite programming 
relaxations of several NP-hard 
graph optimization problems. 
Solving these relaxations is a computational challenge because 
of the potentially large number of violated subgraph constraints.
We introduce a computational framework for these relaxations designed
to cope with these difficulties. We suggest a partial 
Lagrangian dual, and exploit the fact that its
evaluation  decomposes into two independent subproblems.  
This opens the way to use the bundle method from non-smooth 
optimization to minimize the dual function.
Computational experiments on the Max-Cut, stable set and coloring
problem show the efficiency of this approach.

\end{abstract}
%
%
%

%%%%%%%%%%%%%%%%%%%%%%%%%%%%%%%%%%%%%%%%%%%%%%%%%%%%%%%%%%%%%%%%%%%%%%%%%%%%%%%%%%%%%%%%%%%%%%%%%%%%%%%%%%%%%%%%%%%%%%%%%%%%%%%%%%%%%%%%%%%%%%%%%%%%%%%%%%%%%%%%%%%%

\section{Introduction}
The study of NP-hard problems has led to the introduction of 
various hierarchies of relaxations, which typically involve several 
levels. Moving from one level to the next the relaxations get 
increasingly tighter and ultimately the exact optimum may be reached,
but the computational effort grows accordingly.

Among the most prominent hierarchies are the polyhedral ones from Boros, Crama and Hammer \cite{BorosCramaHammer} as well as the ones from Sherali and Adams \cite{SheraliAdamsHierarchy}, Lov\'{a}sz and Schrijver \cite{LovaszSchrijverHierarchy} and Lasserre \cite{LasserreHierarchy} which are based on semidefinite programming (SDP). Even though on the starting level they have a simple SDP relaxation, 
already the first nontrivial level in the hierarchy requires the 
solution of SDPs in matrices of order $\binom{n}{2}$ and on level $k$ the matrix order is $n^{O(k)}$. 
Hence they are considered mainly as theoretical tools and from a practical point of view  these hierarchies 
are of  limited use.

Not all hierarchies are of this type. In \cite{BorosCramaHammer}, a polyhedral
hierarchy for the Max-Cut problem is introduced which maintains 
$\binom{n}{2}$ variables in all levels, with a growing number of 
constraints. 
%Similarly Bonato et al \cite{Bonato} tighten 
%polyhedral Max-Cut relaxations with cutting planes derived from
%small subgraphs. 
More recently, 
Adams, Anjos, Rendl 
and Wiegele~\cite{AARW} introduced a hierarchy of SDP relaxations
which act in the space of symmetric $n \times n$ matrices and at level 
$k$ of the hierarchy all submatrices of order $k$ have to be 'exact' in a 
well-defined sense, i.e. they have to fulfill an \emph{exact subgraph constraint} (ESC).

It is the main purpose of this paper to describe an efficient way
to optimize over level $k$ of this hierarchy for small values of 
$k$, e.g. $k\leqslant 6$, and demonstrate the efficiency
of our approach for the Max-Cut, stable set
and coloring problem. 

Maintaining $\binom{n}{k}$ possible ESCs in 
an SDP in matrices of order $n$ is computationally 
infeasible even for $k=2$ or $k=3$, because each ESC
creates roughly $\binom{k}{2}$ additional equality constraints and at most $2^k$ additional linear variables.

We suggest
the following ideas to overcome this difficulty. First we proceed
iteratively, and in each iteration we include only
(a few hundred of) the most  violated ESCs. 
More importantly, we propose to solve the dual of the resulting
SDP. The structure of this SDP with ESCs
admits a reformulation of the dual in the form of a non-smooth
convex minimization problem with attractive features.
First, any dual solution
yields a valid bound for our relaxations, so it is not necessary 
to carry out the minimization to optimality.
Secondly, the dual function 
evaluation decomposes into two independent problems. The first one
is simply a sum of max-terms (one for each subgraph constraint), 
and the second one consists in solving a 'basic' SDP, independent
of the ESCs. The optimizer for this second problem
also yields a subgradient of the objective function. With this 
information at hand we suggest to use the bundle method from non-smooth 
convex optimization. It 
provides an effective machinery to get close to a minimizer
in few iterations. 

As a result we are able to get near optimal solutions where 
all ESCs for  small values of $k$ ($k \leqslant 6$)
are satisfied up to a small error tolerance. Our computational
results demonstrate the practical potential of this approach.

We finish this introductory section with some notation. 
We denote the vector of all-ones of size $n$ with $\allones{n}$ and $\simplex_n = \{x \in \R^{n}_{+}: \sum_{i=1}^{n}x_{i} = 1\}$. If the dimension is clear from the context we may omit the index and write $\allones{}$ and $\simplex$. Furthermore let $N = \{1, 2, \dots, n\}$. 
A graph $G$ on $n$ vertices has vertex set $N$ and edge set $E$ and $\overline{G}$ is its complement graph. $\Sym{n}$ is the set of $n$-dimensional symmetric matrices. 
%For a vector $x$ we denote with $x_i$ the $i$-th component of $x$. If we have an indexed vector $x_j$, we denote the $i$-th component with $[x_j]_i$.

%%%%%%%%%%%%%%%%%%%%%%%%%%%%%%%%%%%%%%%%%%%%%%%%%%%%%%%%%%%%%%%%%%%%%%%%%%%%%%%%%%%%%%%%%%%%%%%%%%%%%%%%%%%%%%%%%%%%%%%%%%%%%%%%%%%%%%%%%%%%%%%%%%%%%%%%%%%%%%%%%%%%

\section{The Problems and their Semidefinite Relaxations}
\label{sec:DefProblemsBasicRel} 

%In this section we will define the NP-hard Max-Cut, stable set and coloring pro%blem as well as their basic semidefinite relaxations.

%\subsection{The Max-Cut Problem}
In the Max-Cut problem a symmetric matrix $L \in
\Sym{n}$ is given and $c \in \{-1,1 \}^{n}$ which 
maximizes $c^{\transposedVec}Lc$ should be determined. If the matrix $L$ corresponds to the 
Laplacian matrix of a (edge-weighted undirected) graph $G$, 
this is equivalent to finding a bisection of the vertices of $G$ such
that  the total weight of the edges joining the two bisection blocks
is maximized. Such an edge set is also called  a \emph{cut} in  $G$.

Bisections of $N$ can be expressed as $c \in \{-1,1 \}^{n}$ where the 
two bisection blocks correspond to the entries in $c$ of the same sign. 
Given $c \in \{-1,1\}^{n}$ we call $C=cc^{\opvt}$ a \emph{cut
    matrix}. 
The convex hull of all cut matrices (of order $n$) is denoted by
$\CUT_{n}$ or simply $\CUT$ if the dimension is clear. 
Since $c^{\opvt}Lc = \langle L, cc^{\opvt} \rangle$  
Max-Cut can also be written as the following (intractable) linear
program 
$$
z_{mc} = \max \{ \langle L, X\rangle:~ X \in \CUT \}.
$$
$\CUT$ is contained in the spectrahedron 
$
    \setXmc = \left\{ X \in \Sym{n} : \diag(X) = 
    \allones{n},\, X \succcurlyeq 0 \right\},
$
hence
\begin{equation}
\label{relaxation mc}
\max \left\{\langle L,X\rangle :~ X \in \setXmc \right\}
\end{equation}
is a basic semidefinite relaxation for Max-Cut. 
This model is well-known, attributed to Schrijver and  was introduced in a dual 
form by Delorme and Poljak~\cite{DelormePoljak}. 
It can be solved in polynomial time to a fixed prescribed precision and 
solving this relaxation for $n=1000$ takes only a few seconds.

It is well-known that the Max-Cut problem is NP-hard. 
On the positive side, Goemans and Williamson \cite{GoemansWilliamson}
show that one can find  a cut in a graph with
nonnegative edge weights of value at least 0.878$z_{mc}$ in polynomial time. 
%In our setting we allow arbitrary edge weights, and therefore
%consider  arbitrary
%symmetric matrices $L$ as input for the problem.  

%%\begin{equation}\label{conv hull mc}
%%X \in \conv \{C: C \mbox{ cut matrix} \}.
%%\end{equation}

%\subsection{The Stable Set Problem}

In the stable set problem the input
is an unweighted graph $G$. 
We call a set of vertices \emph{stable}, if no two vertices are adjacent.
Moreover we call a vector $s \in \{0,1 \}^n$ a \emph{stable set vector} if it 
is the incidence vector of a stable set. 
The convex hull of all stable set
vectors of $G$  is denoted with $\STAB(G)$.  
In the stable set problem we want to determine the \emph{stability number $\alpha(G)$}, which denotes the cardinality of a largest stable set in $G$, hence
$
\alpha(G) = \max\left\{ \allones{}^{\opvt}s:~ s \in \STAB(G) \right\}.  
$
Furthermore we denote with
%\begin{equation*}
$\STAB^{2}(G) =
\conv \left\{ ss^{\opvt}:~ s \in \STAB(G) \right\}$ 
%\end{equation*}
the convex hull of all \emph{stable set matrices} $ss^{\opvt}$. Then with the arguments of Gaar~\cite{elli-diss} it is easy to check that
$
\alpha(G) = \max\{ \trace(X):~ X \in \STAB^2(G) \}.
$
%Now in order to get a relaxation let $S = \sum_i \alpha_i s_i s_i^\opvt \in \STAB^{2}(G)$ for some stable set vectors $s_i$ and some $\alpha_i \geqslant 0$ with $\sum_i \alpha_i = 1$.
%Then the following properties can easily be shown to hold.
%\begin{equation}\label{def-stab}
%\diag(S) = s,~ 
%S_{ij}= 0~ \forall \{i,j\} \in E,~ 
%\sum_i \alpha_i {1 \choose s_i }{1 \choose s_i}^{\opvt} = 
%\left(
%\begin{array}{cc}
%1 & s^{\opvt} \\
%s & S
%\end{array}
%\right)
%\succcurlyeq 0
%\end{equation}
%So clearly
Furthermore
$\STAB^{2}(G)$ is contained  
in the following spectrahedron
\begin{equation*}
    \setXss = \left\{ X \in \Sym{n} :~ 
    X_{ij}=0 \quad 
    \forall \{i,j\} \in E,~ 
    x = \diag(X),~ 
    \left(
    \begin{array}{cc}
        1 & x^{\opvt} \\
        x & X
    \end{array}
    \right)
    \succcurlyeq 0 \right\},
\end{equation*}
which is known as the \emph{theta body} in the literature. Therefore 
\begin{equation}
\label{relaxation ss}
\vartheta(G)= \max \left\{ \trace(X):~ X \in \setXss \right\}
\end{equation}
is a relaxation of the stable set problem. The Lov{\'a}sz theta function $\vartheta(G)$ was introduced in a seminal paper by Lov{\'a}sz \cite{LovaszStart}.
We refer to Gr{\"o}tschel, Lov{\'a}sz and Schrijver 
\cite{OurUsedFormOfLovasTheta} for a comprehensive analysis of $\vartheta(G)$. 

Determining $\alpha(G)$ is again NP-hard. Contrary to Max-Cut, which
has a polynomial time .878-approximation, for every $\varepsilon>0$ there can be no polynomial time algorithm that 
approximates $\alpha(G)$ within a factor better than $O(n^{1-\varepsilon})$ unless
$P=NP$, see H{\aa}stad \cite{stableSetNotApproximable}.

%\subsection{The Vertex Coloring Problem}

The coloring problem for a given graph $G$ consists in determining 
the \emph{chromatic number} $\chi(G)$, which is the smallest $t$ such that
$N$ can be partitioned into $t$ stable sets. 
Let $S=(s_{1}, \ldots, s_{k})$ be a matrix
where each column is a stable set vector 
and these stable sets 
partition
$V$ into $k$ sets.  
Let us call such  matrices $S$ {\em stable-set partition matrices} (SSPM). 
The $n \times n$ matrix $X=SS^{T}$ is called {\em coloring matrix}. 
%It is a direct sum of $k$ all-ones blocks and each block (=color class)
%corresponds to  a stable set in $G$. 
The convex hull of the
set of all coloring matrices of $G$
is denoted by $\COL(G)$.
%$\COL(G) = \conv \left\{X:~ X \text{ coloring matrix of }G \right\}.$
%As the situation is more involved for the coloring problem, 
We also need the \emph{extended coloring polytope}
\begin{equation*}
\COLE(G) = \conv \left\{
\left( 
\begin{array}{cc} k & \allones{}^{\opvt}\\\allones{} & X
\end{array}\right)
= 
\sum_{i=1}^{k} 
\binom{1}{s_{i}}
\binom{1}{s_{i}}^{\opvt} 
:
\begin{array}{c}
S = (s_{1}, \ldots, s_{k}) \text{ is a} \\
\text{SSPM of } G,~ X = SS^{\opvt} 
\end{array}
\right\}.
\end{equation*}

%As before, we extend each 
%$s_{i}$ by an additional entry which we set to one and consider the
%positive semidefinite matrix
%\begin{equation*}
%\sum_{i} 
% {1 \choose s_{i}}
% {1 \choose s_{i}}^{\opvt}
%=
%\left( 
%\begin{array}{cc} k & \allones{}^{\opvt}\\\allones{} & X
%\end{array}\right) \succcurlyeq 0.
%\end{equation*}
%Thus  
%$\sum_{i}s_{i}= S\allones{}=\allones{}$, because each vertex is in exactly one set. 
%Moreover
%$(SS^{T})_{ij}=0$ for all $\{i,j\} \in E$, since $i$ and $j$ can not be in 
%the same set if there is an edge joining them.  
%Here we used $X=SS^{T}$ and $S\allones{n}=\allones{n}$.
%A matrix $X \in \COL(G)$  satisfies 
%$\diag(X)=\allones{}, X_{ij}=0 ~\forall \{i,j\} \in E$.
%Note also that $rank(S)=k$ implies that $rank(X)=k$ and
%furthermore
%\begin{equation*}
%\left( 
%\begin{array}{cc} t & \allones{}^{\opvt}\\\allones{} & X 
%\end{array}\right) \succcurlyeq 0 ~\Leftrightarrow ~t \geqslant rank(X)=k. 
%\end{equation*}

%Like before in the case of the stable set problem, 
The difficult set $\COLE$ can be relaxed to the easier spectrahedron $\setXcol$
\begin{equation*} %\label{th2c}
\setXcol = \left\{
\left( 
\begin{array}{cc} t & \allones{}^{\opvt}\\\allones{} & X 
\end{array}\right) \succcurlyeq 0:~ \diag(X)=\allones{n}, ~X_{ij}=0 ~\forall \{i,j\} \in E 
\right\}
\end{equation*}
and we can consider the semidefinite program  
\begin{equation}\label{relaxation col}
t^{*}(G) = \min \left\{t:~ 
\left( 
\begin{array}{cc} t & \allones{}^{\opvt}\\\allones{} & X 
\end{array}\right) \in \setXcol \right\}.
\end{equation}
Obviously 
$t^{*}(G) \leqslant \chi(G)$ holds
because the SSPM $S$ consisting of $\chi(G)$ stable sets 
yields a feasible coloring matrix $X=SS^{\opvt}$ with objective function value $\chi(G)$.
It is in fact a consequence of conic duality that
$t^{*}(G)= \vartheta(\overline{G})$ holds.

%This also shows the 'sandwich theorem' of Lov\'asz \cite{LovaszStart}, stating that
%\begin{equation*}
%\alpha(\overline{G}) = \omega(G) \leqslant \vartheta(\overline{G}) \leqslant \chi(G),
%\end{equation*}
%so the tractable parameter $\vartheta(G)$ is 'sandwiched' 
%between the
%NP-hard parameters $\chi(\overline{G})$ and $\alpha(G)$. Here $\omega(G)$ is the clique number of $G$.

It is NP-hard to 
find  $\chi(G)$, to find a 4-coloring of a
3-colorable graph \cite{GuruswamiKhanna} and to color
a $k$-colorable graph with $O(k^{\frac{\log k}{25}})$ colors 
for sufficiently large $k$, \cite{Khot}.

%%%%%%%%%%%%%%%%%%%%%%%%%%%%%%%%%%%%%%%%%%%%%%%%%%%%%%%%%%%%%%%%%%%%%%%%%%%%%%%%%%%%%%%%%%%%%%%%%%%%%%%%%%%%%%%%%%%%%%%%%%%%%%%%%%%%%%%%%%%%%%%%%%%%%%%%%%%%%%%%%%%%

\section{Exact Subgraph Hierarchy}\label{sec: ESC}
In this section we will discuss how to systematically tighten the relaxations \eqref{relaxation mc},
\eqref{relaxation ss} and \eqref{relaxation col} with  'exactness conditions' imposed on small subgraphs. We obtained these relaxations by relaxing the feasible regions $\CUT$, 
$\STAB^{2}$ and $\COL$ of the integer problem to simple spectrahedral sets. Now we will use small subgraphs to get closer to original feasible regions again.

%Towards this end we need the following definitions. 
For $I \subseteq N$ 
we denote with $\XI$ the principal submatrix of $X$ corresponding to the rows and
columns in $I$. 
Furthermore let $\GI$ be the induced subgraph of $G$ 
on the set of vertices $I$ and let $\kI=|I|$ be the cardinality of $I$.

We first look at the exact subgraph relaxations for Max-Cut. 
The \emph{exact subgraph constraint} (ESC) on $I \subseteq N$, introduced
in~\cite{AARW} by Adams, Anjos, Rendl 
and Wiegele,  
requires  
that the matrix $\XI$ corresponding to the subgraph $\GI$ lies in 
the convex hull of the cut matrices of $\GI$, 
that is
\begin{equation*}
%\label{esc mc}
\XI \in  \CUT_{|I|}. 
\end{equation*}
In this case we say that $X$ is \emph{exact} on $I$.
%We recall that there are $\tI=2^{\kI-1}$  cut matrices $C^{I}_{1}$, \dots, $C^{I}_{\tI}$ (obtained through cut vectors $c \in \{-1,1 \}^{\kI}$) for $I$.

Now we want the ESCs to be fulfilled not only 
for one but for a certain selection of  subgraphs.
We denote with $J$ the set of subgraphs which we require to be exact
and  get the following 
SDP relaxation with ESCs for Max-Cut.
\begin{equation}
\label{relaxation mc with esc}
\max \{\langle L,X\rangle:~
X \in \setXmc,~ \XI \in \CUT_{|I|} ~ \forall I \in J \}
\end{equation}

We proceed analogously for the stable set problem in a graph $G$. 
The ESC of a subgraph $\GI$ for 
the stable set problem requires that 
$
\XI \in \STAB^{2}(G_{I})
$
holds and 
%which analogously to Max-Cut makes sure that the matrix 
%$\XI$ corresponding to the subgraph $\GI$ is in the 
%convex hull of all stable set matrices.
%Let $S^{I}_{1}$, \dots, $S^{I}_{\tI}$ denote all stable set matrices 
%of the subgraph $\GI$. 
%Note that the number of stable set matrices $\tI$ 
%of two graphs with the same number of vertices can differ, 
%but $\tI \leqslant 2^{\kI}$ always holds. 
%Analogously to the Max-Cut problem the exact subgraph constraint 
%for the subgraph $\GI$ makes sure that ,  
%\begin{equation*}
%\XI \in  \conv \{ S_{i}^{I}: 1 \leqslant i \leqslant \tI \}.
%\end{equation*}
%Let us again collect all the subgraphs which we require to be exact in the set 
%$J$. 
the SDP with ESCs for the stable set problem is
\begin{equation}
\label{relaxation ss with esc}
\max \{ \trace(X):~ X \in \setXss,~
\XI \in \STAB^{2}(\GI) ~ \forall I \in J \}.
\end{equation}

Turning to  the coloring problem, we analogously impose additional constraints
of the form
$
\XI \in \COL(G_I)
$
to obtain the SDP with ESCs
\begin{equation}
\label{relaxation col with esc}
\min \left\{t:~ 
\left( 
\begin{array}{cc} t & \allones{}^{\opvt}\\ \allones{} & X 
\end{array}\right) \in \setXcol,~
\XI \in \COL(\GI) ~ \forall I \in J
\right\}.
\end{equation}

Note that in the case of the stable set and the coloring problem the polytopes $\STAB^{2}(\GI)$ and $\COL(\GI)$ depend on the subgraph $\GI$, whereas in Max-Cut the polytope $\CUT_{|I|}$ only depends on the number of vertices of the subgraph.

From a theoretical point of view, we obtain the $k$-th level of the exact subgraph hierarchy of \cite{AARW}  if we use $J = \{I \subseteq N:~ |I| = k\}$ in the relaxations~\eqref{relaxation mc with esc}, \eqref{relaxation ss with esc} and~\eqref{relaxation col with esc} respectively. We denote the corresponding objective function values with $\eshmc{k}$, $\eshss{k}$ and $\eshc{k}$. So the $k$-th level of the hierarchy is obtained by forcing all subgraphs on $k$ vertices to be exact in the basic SDP relaxation.

In the case of the stable set and the Max-Cut problem we have $\eshss{n} = \alpha(G)$ (see~\cite{elli-diss}) and $\eshmc{n} = z_{mc}$. For coloring $\eshc{n} \leqslant \chi(G)$ holds. Let $\eshcE{k}$ be the resulting value if we add the inequalities $t \geqslant \sum_{i=1}^{\tI}[\lambdaI]_{i}|S^{I}_{i}|$ where $|S^{I}_{i}|$ is the number of colors used for the SSPM $S^{I}_{i}$ and $\lambdaI \in \simplex_{\tI}$ is a variable for the convex combination for each subgraph $I$ to the SDP for $\eshc{k}$. Then $\eshcE{n} = \chi(G)$ holds. Since the focus of this paper are computational results we are interested only in the computational results we omit the details and further theoretical investigations.

An important feature of this hierarchy is that the size of the matrix variable remains $n$ or $n+1$ on all levels of the hierarchy and only more linear variables and constraints (enforcing the ESCs, hence representing convex hull conditions) are added on higher levels. So it is possible to approximate $\eshmc{k}$, $\eshss{k}$ and $\eshc{k}$ by forcing only some subgraphs of order $k$ to be exact. This is our key ingredient to computationally obtain tight bounds on $z_{mc}$, $\alpha(G)$ and $\chi(G)$.

From a practical point of view solving the relaxations~\eqref{relaxation mc with esc}, \eqref{relaxation ss with esc} 
and~\eqref{relaxation col with esc} with standard interior point (IP) solvers like 
%SeDuMi~\cite{SeDuMi}, 
SDPT3~\cite{SDPT32} or MOSEK~\cite{mosek} is very time consuming. 
In \autoref{MC and SS solution times exact} we list computation times (in seconds) for one specific Max-Cut and one specific stable set instance. 
We vary the number of ESCs for subgraphs of order $3$, $4$ and $5$, so we solve \eqref{relaxation mc with esc} and \eqref{relaxation ss with esc} 
for different $J$. 
We choose $J$ such that the total number of equality constraints induced 
by the convex hull formulation of the ESCs $b$ ranges between 6000 and 15000. 
Since the matrix order $n$ is fixed to $n=100$, the overall computation time
depends essentially on the number of constraints, 
independent of the specific form of the objective function. 
Aside from the ESC constraints, we have $n$ additional equations for Max-Cut
and $n+m+1$ additional equations for the stable set problem. 
Here $m$ denotes the number of edges of the graph. We have $m=722$  in the 
example graph. 
Clearly the running times get huge for a large number of ESC. 
Furthermore MATLAB requires 12 Gigabyte of memory for $b=15000$, 
showing also memory limitations.

Note that it is argued in \cite{AARW} that 
$\eshmc{4}=\eshmc{3}$, so we omit subgraphs of
order $\kI=4$ for Max-Cut. 
This is because in the back of our minds our final algorithm 
to determine the best
possible bounds first includes ESCs of size $k$, 
starting for example with $k=3$. As soon as
we do not find violated ESCs of size $k$ anymore, 
we repeat this for size $k+1$.

%%%%%%%%%%%%%%%%%%%%%%%%%%%%%%%%%%%%%%%%%%%%%%%%%%%%%%%%%%%%%%%%%%%%%%%%%%%%%%%%%%%%%%%%%%%%%%%%%%%%%%%%%%%%%%%%%%%%%%%%%%%%%%%%%%%%%%%%%%%%%%%%%%%%%%%%%%%%%%%%%%%%

\section{Partial Lagrangian Dual} \label{sec: lagrangianDual}

To summarize we are interested in solving relaxations~\eqref{relaxation mc with esc}, \eqref{relaxation ss with esc} and~\eqref{relaxation col with esc} with a potentially large number of ESCs, where using interior point solvers is too time consuming. In this section we will first establish a unified formulation of the relaxations~\eqref{relaxation mc with esc}, \eqref{relaxation ss with esc} and~\eqref{relaxation col with esc}. Then we will build the partial Lagrangian dual of this formulation, where only the ESCs are dualized. This model will be particularly amenable for the bundle method, 
because it will be straightforward to obtain a subgradient 
of the model when evaluating it at a certain point.

In order to unify the notation for the three problems observe that the ESCs $\XI \in \CUT_{|I|}$, $\XI \in \STAB^{2}(\GI)$ and 
$\XI \in \COL(\GI)$ can be represented as
\begin{equation}\label{esc in conv hull}
\XI = \sum_{i=1}^{\tI} \lambda_{i}C^{I}_{i},\quad \lambda \in \simplex_{\tI},
\end{equation}
where $C^{I}_{i}$ is the $i$-th cut, stable set or coloring matrix of the subgraph $\GI$ and $\tI$ is their total number.

A formal description of ESC in \eqref{esc in conv hull} 
requires some additional notation. 
First we introduce the projection $\opP\colon \Sym{n} \mapsto \Sym{\kI}$, 
mapping $X$ to the submatrix $\XI$. 
Second we define a map $\opAmc\colon \Sym{\kI} \mapsto \R^{\tI}$, such that its adjoint map $\opAmct\colon \R^{\tI} \mapsto \Sym{\kI}$ is given by $\opAmct(\lambda)=\sum_{i=1}^{\tI}{\lambda_i C_i^{I}}$ and produces a linear combination of 
the cut, stable set or coloring matrices. 
Thus we can rewrite \eqref{esc in conv hull} as 
\begin{equation}
\label{esc mc withAP}
\opAmct(\lambdaI) - \opP(X) = 0, \quad \lambdaI \in \simplex_{\tI}.
\end{equation}

The left-hand side of the matrix equality is a symmetric matrix, of which some entries (depending on which problem we consider) are zero for sure, so we do not have to include all $\kI \times \kI$ equality constraints into the SDP. Let $\bI$ be the number of equality constraints we have to include.
Note that $\bI = \binom{\kI}{2}$, $\bI = \binom{\kI +1}{2} - m_I$ and $\bI = \binom{\kI}{2} - m_I$ for the Max-Cut, stable set and coloring problem respectively, if $m_I$ denotes the number of edges of $\GI$. This is because in the case of the stable set problem we also have to include equations for the entries of the main diagonal contrary to  Max-Cut and the coloring problem.
 Then we define a linear map $\opMmc\colon \R^{\bI} \mapsto \Sym{\kI}$ such that the adjoint operator $\opMmct\colon \Sym{\kI} \mapsto \R^{\bI}$ extracts the $\bI$ positions, for which we have to include the equality constraints, 
into a vector. 
So eventually we can rephrase \eqref{esc mc withAP}
equivalently as
\begin{align*}
\opMmct(\opAmct(\lambdaI) - \opP(X)) = 0, \quad \lambdaI \in \simplex_{\tI},
\end{align*}
which are $\bI+1$ equalities and $\tI$ inequalities.
In consequence all three relaxations \eqref{relaxation mc with esc},~\eqref{relaxation ss with esc} and~\eqref{relaxation col with esc}
have the generic form
\begin{equation}\label{sdp}
z = \max \{ \langle C, \widehat{X} \rangle:~
\widehat{X} \in \setX,~ \lambdaI \in \simplex_{\tI},~ \opMt(\opAt(\lambdaI) - \opP(X)) = 0  ~ \forall I \in J\},
\end{equation}
where $C$, $\setX$, $\opA$, $\opM$ and $\bI$ 
have to be defined problem specific. Furthermore $\widehat{X} = X$ in the case of Max-Cut and stable set and $\widehat{X} = \left( 
\begin{array}{cc} t & \allones{}^{\opvt}\\ \allones{} & X 
\end{array}\right)$ for coloring, but for the sake of understandability we will just use $X$ in the following.
%The computation times in \autoref{MC and SS solution times exact} contain timings of interior point codes for this relaxation.

The key idea to get a handle on problem \eqref{sdp} is to consider the partial
Lagrangian dual where the ESCs (without the constrains $\lambdaI \in \simplex_{\tI}$) are dualized. We introduce a vector of multipliers $\yI$ of size $\bI$ for each $I$ and collect them in $\yJ = (\yI)_{I\in J}$ and also collect $\lambdaJ = (\lambdaI)_{I\in J}$.
The  Lagrangian function becomes
\begin{equation*}
\Lagrangian(X,\lambdaJ,\yJ) = \langle C, X \rangle + \sum_{I \in J}{\langle \yI, \opMt(\opAt(\lambdaI) - \opP(X)) \rangle}
%\\
%& = \langle C, X \rangle + \sum_{I \in J}{\langle \opM(\yI), \opAt(\lambdaI) - \opP(X) \rangle}\\
%    & = \left\langle C - \sum_{I\in J}{\opPt \opM(\yI)}, X \right\rangle + \sum_{I\in J}{\langle \opA \opM(\yI), \lambdaI \rangle}
\end{equation*} 
and standard duality arguments (Rockafellar~\cite[Corollary 37.3.2]{RockafellarConvexAnalysis}) yield
\begin{equation}\label{lagrangianDual}
z = \min_{\yJ} \max_{\substack{X \in \setX\\ \lambdaI \in \simplex_{\tI}}} \Lagrangian(X,\lambdaJ,\yJ).
\end{equation}
For a fixed set of multipliers $\yJ$ the inner maximization becomes 
\begin{equation*} \max_{\substack{X \in \setX\\ \lambdaI \in \simplex_{\tI}}} \left\langle C - \sum_{I\in J}{\opPt \opM(\yI)}, X \right\rangle + \sum_{I\in J}{\langle \opA \opM(\yI), \lambdaI \rangle}.
\end{equation*}

%Thus $\opD$ can be represented by the $\tI \times \bI $ matrix, whose $i$-th row is given by $\opMmct(C^{I}_i)$ or $\opMsst(S^{I}_i)$ depending on whether we consider the Max-Cut or stable set case respectively.

This maximization is interesting in at least two aspects. 
First, it is separable in the sense that the first term depends only on $X$ and the second one only on the separate $\lambdaI$. 
Moreover, if we denote the linear map $\opA\opM(\yI)\colon \R^{\bI} \mapsto \R^{\tI}$ with $\opD$, the second term has an explicit solution, namely 
\begin{align}\label{solutionMaxTerm}
\max_{\lambdaI \in \simplex_{\tI}}\langle \opD(\yI), \lambdaI \rangle = \max_{1 \leqslant i \leqslant \tI}\left[\opD (\yI)\right]_{i}.
\end{align}

In order to consider the first term in more detail, we define the following function. Let $\sumbI = \sum_{I\in J}\bI$ be the dimension of $\yJ$. Then $h\colon \R^{\sumbI} \to \R$ is defined as
\begin{equation}\label{def:hy}
h(\yJ)= \max_{X \in \setX} \left\langle C - \sum_{I\in J}{\opPt \opM(\yI)}, X \right\rangle = \left\langle C - \sum_{I\in J}{\opPt \opM(\yI)}, X^{\ast} \right\rangle,
\end{equation}
where $X^{\ast}$ is a maximizer over the set $\setX$ for $y$ fixed. 
Note that $h(\yJ)$ is convex but non-smooth, but \eqref{def:hy} shows that
$g_I= -\opM^T \opP(X^{\ast})$ is a subgradient of $h$ with respect to $\yI$.  
%Then the following Lemma holds, see Gaar \cite{elli-diss} for a rigorous proof.
%\begin{lemma} \label{Lemma:subgradient}
%    $h(\yJ)$ is convex. Let $X^\ast$ be a maximizer of $h$ at $\yJ^\ast$. 
%    If $\gI = -\opMt\opP(X^{\ast})$, then $\gJ= (\gI)_{I \in J}$ is a subgradient of $h$ at $\yJ^\ast$.
%\end{lemma}
%\begin{proof}
%    The convexity of $h$ follows from the fact that $h$ is a pointwise maximum of linear functions. 
%    Let $X^{\ast}$ be a maximizer of $h$ at $\yJ^\ast$ over $\setX$. Let $\yJ$ be arbitrary. Then 
%    $$
%    h(\yJ) \geqslant \left\langle C - \sum_I {\opPt \opM(\yI)}, X^{\ast} \right\rangle =
%    h(\yJ^\ast) + \sum_I \langle -\opMt \opP(X^\ast), \yI -\yI^\ast \rangle.
%    %h(\yJ^\ast) + \langle g_J, \yJ -\yJ^\ast \rangle.
%    $$
%    %\begin{equation*}
%    %h(\yJ) = \left\langle C - \sum_{I\in J}{\opPt \opM(\yI)}, X^{\ast} \right\rangle = \langle C , X^{\ast} \rangle - \sum_{I\in J}{\langle \yI, \opMt \opP (X^{\ast})\rangle}
%    %\end{equation*}
%    This shows that $-\opMt\opP(X^{\ast})$ is a subgradient of $h$ at $\yJ^\ast$ with respect to $\yI$.
%\end{proof}
By combining \eqref{solutionMaxTerm} and \eqref{def:hy} we can reformulate the partial Lagrangian dual \eqref{lagrangianDual} to
\begin{equation}\label{primal compact for bundle}
z = \min_{\yJ} \left\{ h(\yJ) + \sum_{I\in J}{\max_{1 \leqslant i \leqslant \tI}\left[\opD (\yI)\right]_i}\right\}.
\end{equation}

The formulation \eqref{primal compact for bundle} of the original relaxations~\eqref{relaxation mc with esc},~\eqref{relaxation ss with esc} and~\eqref{relaxation col with esc}  fits perfectly into the bundle method setting 
described by Frangioni and Gorgone 
in~\cite{EasyBundle}, hence we suggest to approach this problem using the bundle method.

%%%%%%%%%%%%%%%%%%%%%%%%%%%%%%%%%%%%%%%%%%%%%%%%%%%%%%%%%%%%%%%%%%%%%%%%%%%%%%%%%%%%%%%%%%%%%%%%%%%%%%%%%%%%%%%%%%%%%%%%%%%%%%%%%%%%%%%%%%%%%%%%%%%%%%%%%%%%%%%%%%%%

\section{Solving \eqref{primal compact for bundle} with the Bundle Method}
\label{sec:bundle}

The bundle method is an iterative procedure for minimizing a convex non-smooth function and firstly maintains the \textit{current center} $\yCenter$, 
which represents the current estimate to the optimal solution, throughout the iterations.
Secondly it maintains the 
bundle of the form $\Bundle = \{ 
(\yBdlExplanation{1}, \hBdlExplanation{1},\gBdlExplanation{1},\XBdlExplanation{1}), \dots, 
(\yBdlExplanation{\bdlSize}, \hBdlExplanation{\bdlSize},\gBdlExplanation{\bdlSize},\XBdlExplanation{\bdlSize}) \}$.  Here 
$y_1, \ldots, y_r$ are the points which we use to set up our subgradient 
model. Moreover $h_i = h(y_i)$, $g_i$ is a subgradient  of $h$ at $y_i$ 
and $X_i$ is a maximizer of $h$ at $y_i$ as in \eqref{def:hy}. 

At the start we select $y_1=\yCenter =0$ and evaluate $h$ at $\yCenter$, which 
yields the bundle 
$\Bundle =\{(y_1,g_1, h_1,X_1)\}$. A general iteration consists of the two steps determining the new \emph{trial point} and evaluating the \emph{oracle}. For determining a new trial point $\yTrial$ the subgradient information of the bundle $\Bundle$ 
translates into the subgradient model
$ %\begin{align}
\label{subgradientModel}
h(y) \geqslant \hBdlExplanation{j} + \langle \gBdlExplanation{j},y-\yBdlExplanation{j}\rangle$ for all $j = 1$, \dots,  $\bdlSize$.
 %\end{align}
It is common to introduce
$
\eBdlExplanation{j} = h(\yCenter) - \hBdlExplanation{j} - \langle \gBdlExplanation{j},\yCenter-\yBdlExplanation{j}\rangle$ for $j = 1,$ \dots,  $\bdlSize
$
and with $\hCenter = h(\yCenter)$ the subgradient model becomes
\begin{align}
\label{subgrMaxModel}
h(y) \geqslant \max_{1\leqslant j \leqslant r} \left\{\hCenter - \eBdlExplanation{j} + \langle \gBdlExplanation{j},y-\yCenter\rangle \right\}.
\end{align}
The right-hand side above is convex, piecewise linear and minorizes $h$. 
In each iteration of the bundle method we minimize the right-hand 
side of \eqref{subgrMaxModel} instead of $h$, 
but ensure that we do not move too far from $\yCenter$ 
%This is achieved 
by adding a penalty term 
of the form
$
\frac{1}{2}\bdlPar\norm{y-\yCenter}^{2}
$
for a parameter $\bdlPar \in \R_{+}$  
to the objective function.
With the auxiliary variables $w \in \R$ and $\vI \in \R$ for all $I \in J$ to model the maximum terms and with $v = (\vI)_{I \in J} \in \R^\nrESC$ and $\nrESC = |J|$ we  
end up with
\begin{align}
\label{bundleProblemPiO}
\min_{y,w,v} \quad w &+ \sum_{I \in J} \vI + \frac{1}{2}\bdlPar\norm{y-\yCenter}^{2}\\ \nonumber
st \quad w &\geqslant \hCenter - \eBdlExplanation{j} + \langle \gBdlExplanation{j},y-\yCenter\rangle && \forall j = 1,\dots,\bdlSize \\
\nonumber
\vI &\geqslant \left[\opD(\yI)\right]_{i} && \forall i = 1,\dots,\tI \quad \forall I \in J.
\end{align}
This is a convex quadratic problem in $1+\nrESC+\sumbI$ variables 
with $\bdlSize+\sum_{I \in J}\tI$ linear inequality constraints. 
Its solution $(\yTrial,\wTrial,\vTrial)$ includes 
the new trial point $\yTrial$. Problems of this type can be solved efficiently in various ways, see \cite{elli-diss} for further details. In our implementation we view \eqref{bundleProblemPiO} as a rotated second order cone program with one second-order cone constraint and solve it with MOSEK.

%Note that in contrast to the original bundle method, where a subgradient model is used for the whole function, we follow the approach of Frangioni and Gorgone~\cite{EasyBundle} and only have a subgradient model for the 'complicated part' $h(y)$ and keep the second easy sum part of the function we want to minimize explicitly. 

The second step in each bundle iteration is to evaluate the dual function $h$ at $\yTrial$. In our case determining $h(\yTrial)$ means solving the basic SDP relaxation as introduced in Section \ref{sec:DefProblemsBasicRel} with a modified objective function. Hence in the case of Max-Cut the oracle can be evaluated very quickly, whereas evaluating the oracle is computationally more expensive for the stable set and the coloring problem.

The bundle iteration finishes by deciding whether $\yTrial$ becomes the 
    new center (serious step, roughly speaking if the increase of the objective function is good) or not (null step). 
    In either case the new point is included in the bundle, 
    some other elements of the bundle are possibly removed, the bundle parameter $\bdlPar$ is updated  
    and a new iteration starts.

%%%%%%%%%%%%%%%%%%%%%%%%%%%%%%%%%%%%%%%%%%%%%%%%%%%%%%%%%%%%%%%%%%%%%%%%%%%%%%%%%%%%%%%%%%%%%%%%%%%%%%%%%%%%%%%%%%%%%%%%%%%%%%%%%%%%%%%%%%%%%%%%%%%%%%%%%%%%%%%%%%%%

\section{Computational Results and Conclusions}
We close with a small sample of computational results and start with comparing our bundle method with interior point methods.
In our context we are mostly interested to improve
the upper bounds quickly, so we do not run the bundle method described in Section~\ref{sec:bundle} until
we reach a minimizer, but stop after a fixed number of iterations, say $30$. In Table \ref{MC and SS solution times exact} one sees that the running times decrease drastically if we use the bundle method. For $b\approx 15000$ it takes the bundle method only around $8\%$ of the MOSEK running time to get as close as $95\%$ to the optimal value, which is sufficient for our purposes. One sees that our bundle method scales much better for increasing $|J|$.

If we are given a graph and want to get an approximation on $\eshmc{k}$, $\eshss{k}$ and $\eshc{k}$, then we iteratively perform a fixed number, say 30, iterations of the bundle method and then update the set $J$.
We denote the exact subgraph bounds (ESB) obtained in this way with $\esbmc{k}$, $\esbss{k}$ and $\esbc{k}$.

For the sake of brevity we will only outline how to determine $J$ heuristically, see~\cite{elli-diss} for details.
Let $X^\ast$ be the current solution of~\eqref{relaxation mc with esc}, \eqref{relaxation ss with esc} or~\eqref{relaxation col with esc}. We use the fact that the inner product of $\XI^\ast$ and particular matrices of size $\kI$ is potentially small whenever $\XI^\ast$ is not in $\STAB^2(\GI)$. Minimizing this inner product over all subgraphs of order $\kI$ would yield a quadratic assignment problem, so we repeatedly use a local search heuristic for fixed particular matrices in order to obtain potential subgraphs. Then we calculate the projection distances from $\XI^\ast$ to $\STAB^2(\GI)$ for all these subgraphs and include those in $J$ which have the largest distances and hence are violated most.

Finally we present several computational results for obtained ESBs. Note that we refrain 
from comparing the running times of our bundle method with the running time of inter point 
methods, because interior point methods would reach their limit very soon. Hence the bounds 
presented can only be obtained with our methods in reasonable time.

When considering Max-Cut the graphs in Table \ref{MC table 1} are from the Biq Mac library \cite{BiqMacHomepage} 
with $n=100$ vertices. The edge density is 10\%, 50\% and 90\%.
The first 3 instances have positive weights and the remaining 3 have also some negative weights. 
The column labeled $3$ provides the deviation 
(in \%) of  
the ESB with $k=3$ from $z_{mc}$. 
Thus if $p$ is the value in the column
labeled $3$, then 
$\esbmc{3} = (1 + p/100)z_{mc}$. 
The columns labeled 5 and 7 are to be understood in a similar way
for $k=5$ and $k=7$. We note that the improvement of the bound  
from column 3 to column 7 is
quite substantial in all cases. 
We also point out that the relative gap is much larger if also negative
edge weights are present.

In Table \ref{MC table 2} we look at graphs from the Beasley 
collection \cite{BiqMacHomepage} with 
$n=250$. These instances were used by Rendl, Rinaldi and Wiegele \cite{RendlRinaldiWiegele} 
in a Branch-and-Bound setting. 
We only consider the 'hardest' instances from \cite{RendlRinaldiWiegele}
where the Branch-and-Bound tree has more than 200 nodes. 
The table
provides the gap at the root node 
and also the number of nodes in the Branch-and-Bound tree as
reported in \cite{RendlRinaldiWiegele}.  
The column 7-gap contains the gap  after solving our new relaxation with
ESCs up to size $k=7$. 
We find it remarkable that the first 
instance is solved to 
optimality and the gap in the second instance is reduced by 75 \%
compared to the original gap. 
This implies that using our ESBs would expectedly reduce the very high number of required Branch-And-Bound nodes tremendously.

We conclude that for Max-Cut 
our ESB constitute a substantial
improvement compared to the previously used strongest bounds based
on SDP with triangle inequalities. These correspond to the column 
3-gap.

For the calculations for the stable set and the coloring problem
all instances are chosen in such a way that
$\vartheta(G)$ does not coincide and is not very close to $\alpha(G)$ and $\chi(G)$
respectively.

The instances for  the stable set problem  
are taken partly from the DIMACS challenge \cite{DIMACS1992} with some additional instances from \cite{elli-diss} with $n$ ranging
from 26 to 200. 
Table \ref{tab: ttf ss} contains the new bounds. 
Here the starting point is 
the relaxation  $\vartheta(G)$. 
We carry out
10 cycles of adding ESCs. In each cycle we add at most
200 ESCs, so in the final round we have no more than
2000 ESCs. The column heading indicates 
the order of the subgraphs. Here the improvement of the bounds is smaller 
than in the Max-Cut case, but we see that including larger subgraphs
leads to much tighter bounds. In Table \ref{tab: proj dist ss} we show that our approach also
reduces the largest found projection distance over all subgraphs $\GI$ of $\XI$ to the corresponding $\STAB^2(\GI)$ in the course of the cycles. 
This indicates that the violation of the subgraphs decreases over the cycles and less 
and less subgraphs do not fulfill the ESCs. For example the value $0.000$ for the graph spin5
for $\esbss{2}$ at the end of the cycles means that we did not find a violated subgraph of
order $2$ anymore.

Results for a selection of 
coloring instances from \cite{COLInst} are provided  in Tables \ref{tab: ttf col} and \ref{tab: proj dist col}. 
As in the stable set case there is 
only little improvement using small subgraphs ($k=2$ or 3). The 
inclusion of larger subgraphs ($k=6$) shows the potential of the
exact subgraph approach.

Summarizing, we offer the following conclusions from these preliminary
computational results.

$\bullet$ 
Our computational approach based on the partial Lagrangian dual is very efficient in handling also a large number of ESCs. The dual function 
evaluation separates the SDP part from the ESCs and 
therefore opens the way for large-scale computations. The minimization 
of the dual function is carried out as a convex quadratic optimization problem without any SDP constraints, and therefore is also suitable for a large number of 
ESCs. 

$\bullet$
On the practical side we consider the small ESCs 
for Max-Cut a promising new way to tighten bounds for this problem. 
It will be a promising new project to explore these bounds also in a 
Branch-and-Bound setting.

$\bullet$
Our computational results for 
stable set and coloring confirm the theoretical hardness results
for these problems. Here the improvement of the relaxations is small
for $k\leqslant 3$ but including larger subgraphs yields a noticeable improvement of the bounds. It will be a challenge to extend our approach to larger subgraphs.


 \bibliographystyle{splncs04}
 \bibliography{papers}

\newpage
\appendix 
\section{Tables}

%MATLAB file that starts all the computations that I did
% /home/users/egaar/Documents/MathDissertation/MATLAB/MyCode/tightenThetaFunction/Bundle/SourceCode/BundleESC/SolveExactly/PaperBundleESC/PaperBundleESC_FinalVeryShort/letBothRun.m

%Folder for the IP and Bundle solutions
%/home/users/egaar/Documents/MathDissertation/MATLAB/MyCode/tightenThetaFunction/Bundle/SourceCode/BundleESC/SolveExactly/PaperBundleESC/PaperBundleESC_FinalVeryShort/SummaryBundle

\begin{table}
    \setlength\tabcolsep{4.5pt}
    \centering
    \caption{The running times for one Max-Cut and one stable set instance with different fixed sets of ESCs. The graphs of order $n=100$ are from the Erd\H{o}s-R\'enyi model.}
    \begin{tabular}{|c|rrr|r|rr|rr|rr|}
        \hline
        & & & & & \multicolumn{2}{c|}{interior point} & \multicolumn{4}{c|}{our bundle}
        \\ \hline        
        & \multicolumn{3}{|c|}{\#ESC of size} &$b$ &\multicolumn{2}{c|}{time (sec)}
         & \multicolumn{2}{c|}{time (sec)}   & \multicolumn{2}{c|}{$\%$ of MOSEK}
        \\\hline
        & $3$  & $4$  &  $5$   &          & MOSEK   & SDPT3     
        & oracle & overall & time & value\\\hline
        %   0 &    0 &    0 &     0 &  27.2003 &    0 &    1 &    1\\
        \multirow{ 4}{*}{MC}& 2000  & &    0 &  6000 &    18.37 &   49.22 &  1.01 &  6.05 & 32.93 & 97.20\\
        & 2000  & &  300 &  9000 &    55.24 &  134.78 &  1.18 &  9.33 & 16.90 & 95.02\\
        & 4000  & &    0 & 12000 &   104.56 &  289.78 &  1.71 & 11.13 & 10.64 & 93.66\\
        & 3000  & &  600 & 15000 &   184.43 &  525.85 &  1.56 & 14.83 & 8.04  & 94.54\\        
        \hline 
        \multirow{ 4}{*}{SS}& 1050 &    0 &    0 &  5914 &   23.54 &   79.25 &  7.86 & 10.65 & 45.22 & 98.25 \\
        & 1050 &  212 &   63 &  8719 &   50.11 &  174.33 & 10.61 & 16.52 & 32.96 & 97.89\\
        & 2100 &    0 &    0 & 11780 &  126.40 &  388.07 &  7.43 & 12.27 & 9.71 & 93.65 \\
        & 1575 &  318 &  212 & 14653 &  241.29 &  648.83 & 10.79 & 20.21 & 8.38 & 94.44\\
        \hline
    \end{tabular}
    \label{MC and SS solution times exact}
\end{table}
%\end{minipage}

\begin{table}
    \setlength\tabcolsep{5pt}
    \centering
        \caption{The deviation of the ESB to $z_{mc}$ for several Max-Cut instances.}
        \label{MC table 1}
        \begin{tabular}{|l|rrr|r|}
            \hline
            name       &~~~~~~ 3   & ~~~~~  5    & ~~~~~   7  & $z_{mc}$ \\ \hline
            pw01-100.1 & 0.40 & 0.00 & 0.00 &  2060 \\
            pw05-100.1 & 0.90 & 0.51 & 0.39 & 8045 \\
            pw09-100.1 & 0.58 & 0.38 & 0.31 & 13417 \\
            \hline
            w01-100.1  & 0.13 & 0.00 & 0.00 &  719 \\
            w05-100.1  & 3.91 & 1.41 & 0.85 & 1606 \\
            w09-100.1  & 8.06 & 5.66 & 5.09 & 2096 \\ \hline
        \end{tabular}
\end{table}

\begin{table}
    \setlength\tabcolsep{5pt}
    \centering
    %\begin{center}
        \caption{The gap of the ESB  to $z_{mc}$ for two Max-Cut instances.}
        \label{MC table 2}
        \begin{tabular}{|l|rr|r|r|}
            \hline
            name       & BBnodes   &   root gap    &    7-gap  & $z_{mc}$ \\ \hline
            beas-250-6  & 223 &  1.02 & 0.00 & 41014 \\
            beas-250-8  & 4553 & 2.19 & 0.49 & 35726 \\ \hline
        \end{tabular}
    %\end{center}
\end{table}

%folder with computational results:
% /home/users/egaar/Documents/MathDissertation/MATLAB/Results/tightenThetaFunctionStableSet/Run_2018_10_30_IPCO_Paper/Run_2018_10_30_IPCO_Paper_EGpc13

% folder with table summary:
% /home/users/egaar/Documents/MathDissertation/MATLAB/Results/tightenThetaFunctionStableSet/Run_2018_10_30_IPCO_Paper/Run_2018_10_30_IPCO_Paper_EGpc13_Summary

% matlab file to create summary output:
% 

\begin{table}
    \setlength\tabcolsep{4.5pt}
    \centering
    \caption{Tighten $\vartheta(G)$ towards $\alpha(G)$ for several instances for 10 cycles.}
    \begin{tabular}{|l|rr|r|rrrrr|r|}
        \hline 
        name & $n$ & $m$ & $\vartheta(G)$ & $\esbss{2}$~~ & $\esbss{3}$~~ & $\esbss{4}$~~ & $\esbss{5}$~~ & $\esbss{6}$~~ & $\alpha(G)$ \\  
    \hline 
     CubicVT26\_5 &  26 &    39 & 11.82 & 11.82  & 11.00  & 10.98  & 10.54  & 10.46  & 10 \\ 
     Circulant47\_030 &  47 &   282 & 14.30 & 14.30  & 13.61  & 13.21  & 13.24  & 13.14  & 13 \\ 
     G\_50\_0\_5 &  50 &   308 & 13.56 & 13.46  & 13.13  & 12.96  & 12.82  & 12.67  & 12 \\      
     hamming6\_4 &  64 &  1312 & 5.33 & 4.00  & 4.00  & 4.00  & 4.00  & 4.00  & 4 \\ 
     spin5 & 125 &   375 & 55.90 & 55.90  & 50.42  & 50.17  & 50.00  & 50.00  & 50 \\ 
     keller4 & 171 &  5100 & 14.01 & 13.70  & 13.54  & 13.50  & 13.49  & 13.49  & 11 \\ 
     sanr200\_0\_9 & 200 &  2037 & 49.27 & 49.04  & 48.94  & 48.86  & 48.78  & 48.75  & 42 \\ 
     c\_fat200\_5 & 200 & 11427 & 60.35 & 60.34  & 58.00  & 58.00  & 58.00  & 58.00  & 58 \\ 
    \hline 
    \end{tabular}
    \label{tab: ttf ss}   
\end{table}

\begin{table}
    \setlength\tabcolsep{5pt}
    \centering
    \caption{Maximum found projection distance of $\XI$ to $\STAB^{2}(\GI)$ for the computations of \autoref{tab: ttf ss}.}
    \begin{tabular}{|l|r|rrr|rrr|}
        \hline 
        && \multicolumn{3}{c|}{beginning} & \multicolumn{3}{c|}{end} \\
        \hline
        name & $n$ & 
        $\esbc{2} $ & $\esbc{4}$ & $\esbc{6}$ & 
        $\esbc{2} $ & $\esbc{4}$ & $\esbc{6}$ \\ 
        \hline 
     CubicVT26\_5 &  26  & 0.000 & 0.102  & 0.193  & 0.000   & 0.029   & 0.013\\ 
     G\_50\_0\_5 &  50   & 0.087 & 0.093  & 0.118  & 0.000   & 0.013   & 0.024 \\     
     spin5 & 125         & 0.000 & 0.084  & 0.269  & 0.000   & 0.046   & 0.006   \\ 
     sanr200\_0\_9 & 200 & 0.044 & 0.062  & 0.107  & 0.072   & 0.028   & 0.020 \\  
        \hline 
    \end{tabular}
    \label{tab: proj dist ss}
\end{table}

\begin{table}
    \setlength\tabcolsep{4.5pt}
    \centering
    \caption{Tighten $\vartheta(G)$ towards $\chi(G)$ for several instances for 10 cycles.}
    \begin{tabular}{|l|rr|r|rrrrr|r|}
        \hline 
        name & $n$ & $m$ & $\vartheta(G)$ & $\esbc{2}$ & $\esbc{3}$ & $\esbc{4}$ & $\esbc{5}$ & $\esbc{6}$ & $\chi(G) \leqslant$ \\ 
        \hline 
        myciel4 &  23 &    71 & 2.53 & 2.53  & 2.90  & 2.91  & 3.28  & 3.29  & 5 \\ 
        myciel5 &  47 &   236 & 2.64 & 2.64  & 3.05  & 3.09  & 3.45  & 3.45  & 6 \\ 
        mug88\_1 &  88 &   146 & 3.00 & 3.00  & 3.00  & 3.00  & 3.00  & 3.00  & 4 \\ 
        1\_FullIns\_4 &  93 &   593 & 3.12 & 3.12  & 3.25  & 3.37  & 3.80  & 3.80  & 5 \\      
        myciel6 &  95 &   755 & 2.73 & 2.73  & 3.02  & 3.09  & 3.57  & 3.51  & 7 \\               
        myciel7 & 191 &  2360 & 2.82 & 2.82  & 3.02  & 3.08  & 3.63  & 3.50  & 8 \\       
        2\_FullIns\_4 & 212 &  1621 & 4.06 & 4.06  & 4.32  & 4.38  & 4.66  & 4.68  & 6 \\         
        flat300\_26\_0 & 300 & 21633 & 16.99 & 17.04  & 17.12  & 17.10  & 17.12  & 17.12  & 26 \\ 
        \hline 
        \end{tabular}
        \label{tab: ttf col}
\end{table}

\begin{table}
    \setlength\tabcolsep{5pt}
    \centering
    \caption{Maximum found projection distance of $\XI$ to $\COL(\GI)$ for  the computations of \autoref{tab: ttf col}.}
    \begin{tabular}{|l|r|rrr|rrr|}
        \hline 
        && \multicolumn{3}{c|}{beginning} & \multicolumn{3}{c|}{end} \\
        \hline
        name & $n$ & 
        $\esbc{2} $ & $\esbc{4}$ & $\esbc{6}$ & 
        $\esbc{2} $ & $\esbc{4}$ & $\esbc{6}$ \\ 
        \hline 
        myciel4 &         23 & 0.000  & 0.365 & 0.760  & 0.000  & 0.000 & 0.000 \\ 
        1\_FullIns\_4 &   93 & 0.009  & 0.349 & 0.629  & 0.000  & 0.158 & 0.203 \\ 
        myciel7 &        191 & 0.000  & 0.356 & 0.621  & 0.000  & 0.207 & 0.272 \\ 
        flat300\_26\_0 & 300 & 0.127  & 0.279 & 0.360  & 0.143  & 0.142 & 0.091 \\ 
        \hline 
    \end{tabular}
    \label{tab: proj dist col}
\end{table}

\end{document}